\documentclass[11pt]{article}
\usepackage{mathrsfs}
\usepackage{epic,eepic,epsf,epsfig}
\usepackage{amsfonts,srcltx,mathrsfs}
\usepackage{amsmath}
\usepackage{amssymb}
\usepackage{amsbsy}
\usepackage{graphicx}
\usepackage{amsfonts}
\usepackage{color}
\usepackage{mathrsfs}
\usepackage{epic,eepic,epsf,epsfig}
\usepackage{graphicx}
\usepackage{amsfonts,srcltx,mathrsfs}
\usepackage{array}
\usepackage{amsthm}
\newtheorem{theorem}{Theorem}
\newtheorem{lem}{Lemma}

\newtheorem{defi}{Definition}

\def\f{\noindent}

\begin{document}

\markboth{ et. al}{The $g$-good neighbour diagnosability of hierarchical cubic networks}

\title{The $g$-good neighbour diagnosability of hierarchical cubic networks
}

\author{Shu-Li Zhao, Rong-Xia Hao\footnote{Corresponding author. Email: rxhao@bjtu.edu.cn (R.-X. Hao)}\\[0.2cm]
{\em\small Department of Mathematics, Beijing Jiaotong University,}\\ {\small\em
Beijing 100044, P.R. China}}

\date{}
\maketitle

Let $G=(V, E)$ be a connected graph, a subset $S\subseteq V(G)$ is called an $R^{g}$-vertex-cut of $G$ if $G-F$ is disconnected and any vertex in $G-F$ has at least $g$ neighbours in $G-F$. The $R^{g}$-vertex-connectivity is the size of the minimum $R^{g}$-vertex-cut and denoted by $\kappa^{g}(G)$. Many large-scale multiprocessor or multi-computer systems take interconnection networks as underlying topologies. Fault diagnosis is especially important to identify fault tolerability
of such systems. The $g$-good-neighbor diagnosability such that every fault-free node has at least $g$ fault-free neighbors is a novel measure of diagnosability. In this paper,
we show that the $g$-good-neighbor diagnosability of the hierarchical cubic networks $HCN_{n}$ under the
PMC model for $1\leq g\leq n-1$ and the $MM^{*}$ model for $1\leq g\leq n-1$ is $2^{g}(n+2-g)-1$, respectively.

\medskip

\f {\em Keywords:} Diagnosability; Fault-tolerance; $PMC$ model; $MM^{*}$ model; $g$-good-neighbour diagnosability

\section{Introduction}

With the development of technology, the high performance for large multiprocessor systems is of great importance. For multiprocessor systems, they usually take interconnection networks as underlying
topology. An interconnection network is usually modelled by a connected graph $G=(V, E)$, where vertices represent processors
and edges represent communication links between processors. The connectivity $\kappa(G)$ of a graph $G$ is defined as the minimum number of vertices whose removal disconnects the graph $G$ and the edge connectivity $\lambda(G)$ is defined as the minimum number of edges whose deletion disconnects the graph $G$. They are two important parameters to evaluate the reliability of a network. It was shown in~\cite{X} that the higher these parameters are, the reliable the network is. However, these parameters always underestimate the resilience of a network. To overcome the shortcoming, Esfahanian~\cite{E} introduced the concept of restricted connectivity, which was a parameter to evaluate the fault tolerance of the network in terms of vertex failure. Later, Latifi $et$ $al.$~\cite{La} and Oh and Choi~\cite{O} generalized the parameter to $R^{g}$-vertex-connectivity.

For a connected graph $G=(V, E)$, a subset $S\subseteq V(G)$ is called an $R^{g}$-vertex-cut of $G$ if $G-F$ is disconnected and any vertex in $G-F$ has at least $g$ neighbours in $G-F$. The $R^{g}$-vertex-connectivity is the size of the minimum $R^{g}$-vertex-cut and denoted by $\kappa^{g}(G)$. There are many results about $R^{g}$-vertex-connectivity of networks, one can refer~\cite{Wan, Yang, Yu, Yang1}.

In addition, processors may fail and create faults in the large multiprocessor system. Hence, node fault identification is also of great importance for the system. The first step to deal with faults is to identify the
faulty processors from the fault-free ones. The identification process is called the diagnosis of the system. A system is said to be
$t$-diagnosable if all faulty processors can be identified without replacement, provided that
the number of faults presented does not exceed $t$. The diagnosability $t(G)$ of a system $G$ is the maximum value of $t$
such that $G$ is t-diagnosable~\cite{D,F,L}.

To identify the faulty processors, some diagnosis models were proposed. One of which was introduced by Preparata $et$ $al.$~\cite{Pr} in $1967$ and it is called the $PMC$ diagnosis model. The diagnosis
of the system is achieved through two linked processors testing each other. Another is the $MM^{*}$ diagnosis model, which was proposed by Maeng and Malek~\cite{M} in $1981$. For the $MM^{*}$ model, to diagnose the system, a node sends the same task to two of its neighbours and then compares these responses. In $2005$, Lai $et$ $al.$~\cite{L} introduced the restricted diagnosability of a system, which is called conditional diagnosability. They consider the situation that any fault set can not contain all neighbours of any vertex in the system. In $2012$, Peng $et$ $al.$~\cite{Pe} proposed a new measurement for fault diagnosis of the system, that is, the $g$-good-neighbour diagnosability. This kind of diagnosis requires that every fault-free node contains at least $g$ fault-free neighbours and they studied the $g$-good-neighbour diagnosability of the $n$-dimensional hypercube under the $PMC$ model in~\cite{Pe}. In addition, there are many results about the $g$-good-neighbour diagnosability of other networks. For example, Wang and Lin $et$ $al.$~\cite{w3} studied the $1$-good-neighbour connectivity and diagnosability of Cayley graphs generated by complete graphs; Wang and Han studied the $g$-good-neighbour diagnosability of the $n$-dimensional hypercube under the $MM^{\ast}$ model in~\cite{w}; Yuan $et$ $al.$~\cite{Y} studied the $g$-good-neighbour diagnosability of the $k$-ary $n$-cube under the $PMC$ model and $MM^{*}$ model; Wang and Guo~\cite{w1} studied the $1$-good-neighbour diagnosability of Cayley graphs generated by transposition trees under the $PMC$ model and $MM^{*}$ model; Wang and Lin~\cite{w2} studied the $2$-good-neighbour diagnosability of Cayley graphs generated by transposition trees under the $PMC$ model and $MM^{*}$ model and so on.

In this paper, we focus on the hierarchical cubic networks $HCN_{n}$ for $n\geq 2$. we show that the $g$-good-neighbor diagnosability of the hierarchical cubic networks $HCN_{n}$ under the
PMC model for $1\leq g\leq n-1$ and the comparison model for $1\leq g\leq n-1$ is $2^{g}(n+2-g)-1$, respectively.

\section{Preliminary}
In this section, we will introduce some definitions and notations needed for our discussion. Let $G=(V, E)$ be a non-complete undirected graph, {\em the degree} of a vertex $v\in V(G)$, denoted by $d_{G}(v)$, is the number of edges incident with $v$. The {\em minimum degree} of a vertex $v$ in $G$ is denoted by $\delta(G)$. {\em The neighbourhoods} of the vertex $v$ in $G$ is denoted by $N_{G}(v)$. Let $S\subset V(G),$ we use $N_{G}(S)$ to denote the vertex set $\bigcup_{v\in S}N_{G}(v)\setminus S$. If for any vertex $v\in V(G)$, $d_{G}(v)=k$, then the graph is called {\em $k$-regular}. For any two vertices $u$ and $v$ in $G$, we denote by $cn(G; u, v)$ {\em the number of common neighbours of $u$ and $v$}, that is, $cn(G; u, v)=|N_{G}(u)\bigcap N_{G}(v)|$. Let $F_{1}, F_{2}\subseteq V(G), F_{1}\vartriangle F_{2}=(F_{1}\bigcup F_{2})\setminus (F_{1}\bigcap F_{2})$. The subgraph induced by $V\subseteq V(G)$, denoted by $G[V]$, is a graph whose vertex set is $V$ and the edge set is the set of all the edges of $G$ with both ends in $V$. A fault set $F\subset V(G)$ is called a {\em $g$-good-neighbour faulty set} if for any vertex $v\in V(G)\setminus F$, $|N(v)\bigcap(V\setminus F)|\geq g$.

To diagnose the faults in a system $G$, a number of tests need to perform on the vertices. The collection of all test results is called a syndrome. For
a given syndrome $\sigma$, a subset $F\subset V(G)$ is said to be consistent with $\sigma$ if the syndrome $\sigma$ can
arise from the situation that all vertices in $F$ are faulty and all vertices in $V(G)\setminus F$ are fault-free. If for every syndrome $\sigma$, there is a unique $F\subseteq V(G)$ such that $F$ is consistent with $\sigma$, then the system is said to be diagnosable. Let $\sigma(F)$ denote the set of
all syndromes which are consistent with $F$. Two distinct vertex sets $F_{1}$ and $F_{2}$ in $V(G)$ are said to be indistinguishable if
$\sigma(F_{1})\bigcap \sigma(F_{2})\neq \emptyset $, otherwise, $F_{1}$ and $F_{2}$ are said to be distinguishable. Also, we say $(F_{1},F_{2})$ is an
indistinguishable pair if $\sigma(F_{1})\bigcap \sigma(F_{2})\neq \emptyset$ and a distinguishable pair if $\sigma(F_{1})\bigcap \sigma(F_{2})= \emptyset $.

A system $G=(V,E)$ is $g$-good-neighbour $t$-diagnosable if $F_{1}$ and $F_{2}$ are distinguishable for each
distinct pair of $g$-good-neighbour faulty subsets $F_{1}$ and $F_{2}$ of $V$ with $|F_{1}|\leq t$ and $|F_{2}|\leq t$. The $g$-good-neighbour diagnosability $t_{g}(G)$ of a system $G$ is the maximum value of $t$ such that $G$ is $g$-good-neighbour $t$-diagnosable.

Let $V_{n}$ be the set of binary sequence of length $n$, i.e., $V_{n}=\{x_{1}x_{2}\cdots x_{n}: x_{i}\in\{0,1\}, 1\leq i\leq n\}$. For $x=x_{1}x_{2}\cdots x_{n}\in V_{n}$, the element $\overline{x}=\overline{x}_{1}\overline{x}_{2}\cdots \overline{x}_{n}\in V_{n}$ is called the bitwise complement of $x$, where $\overline{x}_{i}=\{0,1\}\backslash\{x_{i}\}$ for each $i\in\{1,2,\cdots,n\}$.

The hypercube is one of the most fundamental interconnection networks. An $n$-dimensi\\onal hypercube, shortly $n$-cube, is an undirected graph $Q_{n}=(V,E)$ with $|V|=2^{n}$ and $|E|=n2^{n-1}.$ Each vertex can be represented by an $n$-bit binary string. There is an edge between two vertices whenever their binary string representation differs in only one bit position. The hierarchical cubic network was introduced by Ghose and Desai in \cite{gho}, which is feasible to be implemented with thousands of or more processors, with retaining some good properties of the hypercubes, such as regularity, symmetry and logarithmic diameter.  Following, we will introduce the definition of the hierarchical cubic networks.

\begin{defi}\label{defi1} An $n$-dimensional hierarchical cubic network $HCN_{n}$ with vertex set $V_{n}\times V_{n}$ is obtained from $2^{n}$ $n$-cubes $\{xQ_{n}:x\in V_{n}\}$ by adding edges between two $n$-cubes, called crossing edges, according to the following rule: A vertex $(x,y)$ in $xQ_{n}$ is linked to

$(1)$ $(y, x)$ in $yQ_{n}$ if $x\neq y$ or

$(2)$ $(\overline{x}, \overline{y})$ in $\overline{x}Q_{n}$ if $x=y$.

\end{defi}\label{defi1}

The vertex $(y, x)$ in $yQ_{n}$ or $(\overline{x}, \overline{y})$ in $\overline{x}Q_{n}$ is called an outside neighbour of $(x, y)$ in $xQ_{n}$. By the definition of the hierarchical cubic network, it is an $(n+1)$-regular network. Let $x_{1}Q_{n}, x_{2}Q_{n}, \cdots, x_{2^{n}}Q_{n}$ be $n$-cubes of $HCN_{n}$, where $x_{i}\in V_{n}$ for $1\leq i\leq 2^{n}$. An $2$-dimensional hierarchical cubic network $HCN_{2}$ is shown in Fig.1, where the red edges are the crossing edges of $HCN_{2}$.

\begin{figure}[!ht]
\begin{center}
\includegraphics[scale=0.4]{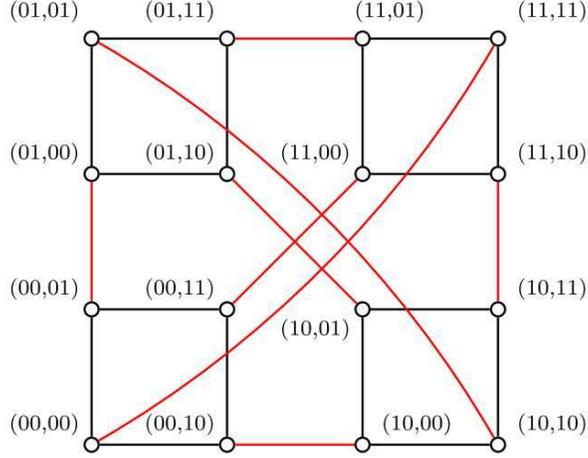}
\end{center}
\caption{The $2$-dimensional hierarchical cubic network $HCN_{2}$}\label{F1}
\end{figure}

There are some results about hierarchical cubic network, one can refer \cite{chi, fu1, fu2, gho, yun1, yun2, zhou} etc. for the detail.

The paper is organized as follows. In section $3$, the $g$-good-neighbour diagnosability of $HCN_{n}$ under the $PMC$ model is determined. In section $4$, the $g$-good-neighbour diagnosability of $HCN_{n}$ under the $MM^{\ast}$ model is determined. In section $5$, the paper is concluded.

\section{The $g$-good neighbour diagnosability of $HCN_{n}$ under the $PMC$ model}
The following properties about the hierarchical cubic networks $HCN_{n}$ are useful.

\begin{lem}\label{lem1} Let $HCN_{n}$ be an $n$-dimensional hierarchical cubic network for $n\geq 2$. Then the following results hold.
\begin{enumerate}
\item [{\rm (1)}] Any vertex of $HCN_{n}$ has exactly one outside neighbour.
\item [{\rm (2)}] There are two crossing edges between two $n$-cubes $xQ_{n}$ and $yQ_{n}$ if and only if $x$ and $y$ are complementary; otherwise there is only one crossing edge.

\item [{\rm (3)}] The set of crossing edges consists of a perfect matching of $HCN_{n}$.
\end{enumerate}
\end{lem}

As $HCN_{n}$ is made up of $2^{n}$ $n$-cubes and a perfect matching, some properties on an $n$-cube $Q_{n}$ are very useful for the proofs of the main results.

\begin{lem}{\rm(\cite{Wu})}\label{lem2}
 If $X$ is a subgraph of $Q_{n}$ and $\delta(X)\geq g$, then $|X|\geq 2^{g}$.
\end{lem}

\begin{lem}\label{lem2.1.1}
Let $HCN_{n}$ be an $n$-dimensional hierarchical cubic network for $n\geq 2$, then there is no triangle in $HCN_{n}$.
\end{lem}
\f {\bf Proof.} Let $x_{1}Q_{n}, x_{2}Q_{n}, \cdots, x_{2^{n}}Q_{n}$ be the $2^{n}$ $n$-cubes of $HCN_{n}$. Suppose to the contrary, that is, there is a triangle in $HCN_{n}$. Let $C=xyzx$ be the triangle in $HCN_{n}$. To prove the result, the following three cases are considered.

Case 1. $x, y, z$ belong to exactly one $n$-cube of $HCN_{n}$.

Without loss of generality, let $x, y, z\in V(x_{1}Q_{n})$. As $Q_{n}$ is bipartite, there is no triangle in $x_{1}Q_{n}$, a contradiction.

Case 2. $x, y, z$ belong to two different $n$-cubes of $HCN_{n}$.

Without loss of generality, let $x, y\in V(x_{1}Q_{n})$ and $z\in V(x_{2}Q_{n})$. As $xz, yz\in E(HCN_{n})$, then $x$ and $y$ are two outside neighbours of $z$, contradict with (1) of Lemma $1$.

Case 3. $x, y, z$ belong to three different $n$-cubes of $HCN_{n}$.

Without loss of generality, let $x\in V(x_{1}Q_{n}), y\in V(x_{2}Q_{n})$ and $z\in V(x_{3}Q_{n})$. As $xz, yz\in E(HCN_{n})$, then $x$ and $y$ are two outside neighbours of $z$, contradict with (1) of Lemma $1$.

Thus, there is no triangle in $HCN_{n}$.

\hfill\qed

\begin{lem}\label{lem2.1}
 If $X$ is a subgraph of $HCN_{n}$ and $\delta(X)\geq g$, then $|X|\geq 2^{g}$ for $1\leq g\leq n-1$.
\end{lem}

\f {\bf Proof.} Let $x_{1}Q_{n}, x_{2}Q_{n}, \cdots, x_{2^{n}}Q_{n}$ be the $2^{n}$ $n$-cubes of $HCN_{n}$. Let $I=\{i|V(X)\cap V(x_{i}Q_{n})\neq\emptyset\}$ for $1\leq i\leq 2^{n}$. To prove the result, the following cases are considered.

Case 1. $|I|=1$

Without loss of generality, let $V(X)\cap V(x_{1}Q_{n})\neq\emptyset$. By Lemma~\ref{lem2}, $|X|\geq 2^{g}$.

Case 2. $|I|\geq 2$

As $V(X)\cap V(x_{i}Q_{n})\neq\emptyset$ for each $i\in [n]$. By the definition of $HCN_{n}$, any vertex of $HCN_{n}$ has exactly one outside neighbour. Thus, $d_{x_{i}Q_{n}}(u)\geq g-1$ for any vertex $u\in V(x_{i}Q_{n})$ and $i\in I$. By Lemma~\ref{lem2}, $|V(X)\cap V(x_{i}Q_{n})|\geq 2^{g-1}$ for each $i\in I$. Thus, $|X|=|\cup_{i\in I} V(X)\cap V(x_{i}Q_{n})|\geq 2^{g-1}+ 2^{g-1}=2^{g}$.
\hfill\qed

\begin{lem}{\rm(\cite{Li})}\label{lem3}
 For $n\geq 1$, $\kappa^{g}(HCN_{n})=2^{g}(n+1-g)$ for any $g$ with $0\leq g\leq n-1$.
\end{lem}

\begin{lem}{\rm(\cite{Li})}\label{lem4}
 Let $x_{1}Q_{n}, x_{2}Q_{n}, \cdots, x_{2^{n}}Q_{n}$ be the $2^{n}$ $n$-cubes of $HCN_{n}$. Let $X\subseteq V(x_{i}Q_{n})$ and $HCN_{n}[X]\cong Q_{g}$. Then $|N_{HCN_{n}}(X)|=2^{g}(n+1-g)$ and $\delta(HCN_{n}[V(HCN_{n})\setminus N_{HCN_{n}}[X]])\geq g$ for $1\leq g\leq n-1$.
\end{lem}

\begin{lem}{\rm(\cite{Zhu})}\label{lem5}
Any two vertices in $V(Q_{n})$ have exactly two common neighbors for $n\geq 3$ if they have any.
\end{lem}

By Lemma~\ref{lem5} and the definition of the hierarchical cubic network $HCN_{n}$, the following lemma holds.

\begin{lem}\label{lem6}
Let $x_{1}Q_{n}, x_{2}Q_{n}, \cdots, x_{2^{n}}Q_{n}$ be the $2^{n}$ $n$-cubes of $HCN_{n}$. Then for any two distinct vertices $u$ and $v$ in the hierarchical cubic network $HCN_{n}$, we have

$$ |N_{HCN_{n}}(u)\cap N_{HCN_{n}}(v)|\left\{
\begin{array}{rcl}
2       &      & if~$u$~and~$v$~belong~to~the~same~$n$~cube~and~ N_{HCN_{n}}(u) \cap\\
&  &   N_{HCN_{n}}(v)\neq\emptyset; \\
1     &      & if~$u$~and~$v$~belong~to~$different$~$n$~cubes~and~ N_{HCN_{n}}(u) \cap\\
&  &   N_{HCN_{n}}(v)\neq\emptyset; \\
0       &      & {otherwise.}
\end{array} \right. $$

\end{lem}

To prove the $g$-good neighbour diagnosability of $HCN_{n}$ under the $PMC$ model, the following results are useful.

\begin{defi}{\rm(\cite{Y})}\label{defi2}
 A system $G=(V,E)$ is $g$-good-neighbour $t$-diagnosable if and only if for any two distinct $g$-good-neighbour faulty subsets $F_{1}$ and $F_{2}$ of $V$ such that $|F_{1}|\leq t$ and $|F_{2}|\leq t$, the sets $F_{1}$ and $F_{2}$ are distinguishable.
\end{defi}

\begin{lem}{\rm(\cite{Y})}\label{lem1.1}
 For any two distinct subsets $F_{1}$ and $F_{2}$ in a system $G=(V,E)$, the sets $F_{1}$ and $F_{2}$ are distinguishable if and only if there exists a vertex $u\in V-(F_{1}\bigcup F_{2})$ and $v\in F_{1}\Delta F_{2}$ such that $(u,v)\in E$.
\end{lem}

\begin{theorem}\label{thm1}
Let $HCN_{n}$ be the $n$-dimensional hierarchical cubic network for $n\geq 2$, then the $g$-good neighbour diagnosability of $HCN_{n}$ under the $PMC$ model satisfies $t_{g}(HCN_{n})\leq 2^{g}(n+2-g)-1$ for $n\geq 2$ and $1\leq g\leq n-1$.
\end{theorem}

\f {\bf Proof.} Let $x_{1}Q_{n}, x_{2}Q_{n}, \cdots, x_{2^{n}}Q_{n}$ be the $2^{n}$ $n$-cubes of $HCN_{n}$. Let $X\subseteq V(x_{1}Q_{n})$ and $HCN_{n}[X]\cong Q_{g}$. Let $F_{1}=N_{HCN_{n}}(X)$ and $F_{2}=N_{HCN_{n}}[X]$. By Lemma~\ref{lem4}, both $F_{1}$ and $F_{2}$ are $g$-good neighbour faulty sets. As $|N_{HCN_{n}}(X)|=2^{g}(n+1-g)$, we have
\vskip0.1cm
~~~~~~$|F_{1}|=2^{g}(n+1-g)\leq 2^{g}(n+2-g)$
\vskip0.1cm
and
\vskip0.1cm
~~~~~~~$|F_{2}|=|F_{1}|+|X|$
\vskip0.1cm
~~~~~~~~~~~~~$=2^{g}(n+1-g)+ 2^{g}$
\vskip0.1cm
~~~~~~~~~~~~~$\leq 2^{g}(n+2-g)$
\vskip0.1cm
As $F_{1}\Delta F_{2}=X$, there is no cross edge between $HCN_{n}\setminus (F_{1}\cup F_{2})$ and $F_{1}\Delta F_{2}$. By Lemma~\ref{lem1.1}, the $g$-good neighbour faulty sets of $F_{1}$ and $F_{2}$ are indistinguishable. By Definition~\ref{defi2}, then hierarchical cubic network $HCN_{n}$ is not $g$-good neighbour $2^{g}(n+2-g)$-diagnosable under the $PMC$ model. Thus, $t_{g}(HCN_{n})\leq 2^{g}(n+2-g)-1$.
\hfill\qed

\begin{theorem}\label{thm2}
Let $HCN_{n}$ be the $n$-dimensional hierarchical cubic network for $n\geq 2$, then the $g$-good neighbour diagnosability of $HCN_{n}$ under the $PMC$ model satisfies $t_{g}(HCN_{n})\geq 2^{g}(n+2-g)-1$ for $n\geq 2$ and $1\leq g\leq n-1$.
\end{theorem}

\f {\bf Proof.} Let $x_{1}Q_{n}, x_{2}Q_{n}, \cdots, x_{2^{n}}Q_{n}$ be the $2^{n}$ $n$-cubes of $HCN_{n}$. To Prove the result, we just need to show that for any two distinct $g$-good neighbour faulty subsets $F_{1}$ and $F_{2}$ of $HCN_{n}$ such that $|F_{1}|\leq 2^{g}(n+2-g)-1$ and $|F_{2}|\leq 2^{g}(n+2-g)-1$, the sets $F_{1}$ and $F_{2}$ are distinguishable. By Lemma~\ref{lem1.1}, there is an edge between $V(HCN_{n})\setminus(F_{1}\cup F_{2})$ and $F_{1}\bigtriangleup F_{2}$.

We prove the result by contradiction. That is, there are two distinct $g$-good neighbour faulty subsets $F_{1}$ and $F_{2}$ of $HCN_{n}$ such that $|F_{1}|\leq 2^{g}(n+2-g)-1$ and $|F_{2}|\leq 2^{g}(n+2-g)-1$, but they are indistinguishable.

First, we show that $V(HCN_{n})\neq F_{1}\cup F_{2}$.

Suppose to the contrary, that is, $V(HCN_{n})=F_{1}\cup F_{2}$. For $1\leq g\leq n-1$, we have

\vskip0.1cm
~~~~~~$|F_{1}\cup F_{2}|=|F_{1}|+|F_{2}|-|F_{1}\cap F_{2}|$

\vskip0.1cm
~~~~~~~~~~~~~~~~~~~$\leq 2[2^{g}(n+2-g)-1]$

\vskip0.1cm
Let $f(g)=2^{g}(n+2-g)$, then $f^{\prime}(g)=2^{g}[(n+2-g)\cdot ln2-1]$. If $f^{\prime}(g)=0$, then $g=n+2-\frac{1}{ln2}> n$. Thus, $f(g)$ is monotonically increasing for $1\leq g\leq n-1$. Thus,

\vskip0.1cm
~~~~~~$f(g)_{max}=f(n-1)$

\vskip0.1cm
~~~~~~~~~~~~~~~~~$=2^{n-1}[n+2-(n-1)]$

\vskip0.1cm
~~~~~~~~~~~~~~~~~$=3\cdot 2^{n-1}$.

Thus, we have

\vskip0.1cm
~~~~~~$|F_{1}\cup F_{2}|\leq 2(3\cdot 2^{n-1}-1)$

\vskip0.1cm
~~~~~~~~~~~~~~~~~$=3\cdot 2^{n}-2$.
\vskip0.1cm
As $|V(HCN_{n})|=2^{2n}$. Obviously, $2^{2n}>3\cdot 2^{n}-2$ for $n\geq 2$, which is a contradiction.

Second, we prove the main result. Without loss of generality, let $F_{2}\setminus F_{1}\neq \emptyset$. As $F_{1}$ is a $g$-good neighbour faulty set, then for any vertex $u$ of $HCN_{n}\setminus F_{1}$, $d_{HCN_{n}\setminus F_{1}}(u)\geq g$. As there is no cross edge between $HCN_{n}\setminus (F_{1}\cup F_{2})$ and $F_{1}\bigtriangleup F_{2}$, then $d_{HCN_{n}[V(HCN_{n})\setminus F_{1}]}(u)=d_{HCN_{n}[F_{2}\setminus F_{1}]}(u)\geq g$. Thus, for any vertex of $u\in F_{2}\setminus F_{1}$, $d_{HCN_{n}[F_{2}\setminus F_{1}]}(u)\geq g$. By Lemma~\ref{lem2.1}, $|F_{2}\setminus F_{1}|\geq 2^{g}$. As $F_{1}$ and $F_{2}$ are both $g$-good neighbour faulty sets, $F_{1}\cap F_{2}$ is also a $g$-good neighbour faulty set. In addition, there is cross edge between $HCN_{n}\setminus (F_{1}\cup F_{2})$ and $F_{1}\bigtriangleup F_{2}$, $F_{1}\cap F_{2}$ is a $g$-good neighbour faulty cut. By Lemma~\ref{lem4}, $|F_{1}\cap F_{2}|\geq 2^{g}(n+1-g)$. Thus

\vskip0.1cm
~~~~~~~~~~~~~$|F_{2}|=|F_{2}\setminus F_{1}|+|F_{1}\cap F_{2}|$

\vskip0.1cm
~~~~~~~~~~~~~~~~~~$\geq 2^{g}+2^{g}(n+1-g)$

\vskip0.1cm
~~~~~~~~~~~~~~~~~~$= 2^{g}(n+2-g)$

By the hypothesis, $|F_{2}|\leq 2^{g}(n+2-g)-1$, which is a contradiction. Thus, $t_{g}(HCN_{n})\geq 2^{g}(n+2-g)-1$ for $n\geq 2$ and $1\leq g\leq n-1$.
\hfill\qed

\vskip0.2cm

By Theorem~\ref{thm1} and Theorem~\ref{thm2}, the following theorem can be obtained.

\begin{theorem}\label{thm2}
Let $HCN_{n}$ be the $n$-dimensional hierarchical cubic network for $n\geq 2$, then the $g$-good neighbour diagnosability under the $PMC$ model is $t_{g}(HCN_{n})= 2^{g}(n+2-g)-1$ for $n\geq 2$ and $1\leq g\leq n-1$.
\end{theorem}

\section{The $g$-good neighbour diagnosability of $HCN_{n}$ under $MM^{*}$ model}

In this section, we will determine the $g$-good neighbour diagnosability of $HCN_{n}$ under the $MM^{*}$ model.

\begin{lem}{\rm(\cite{Y})}\label{lem4.1}
Let $G=(V,E)$ be a system under the $MM^{*}$ model. Two distinct subsets $F_{1}$ and $F_{2}$ of $V$ are distinguishable if and only if     at least one of the following conditions holds:
\begin{enumerate}
\item [{\rm (1)}]\label{4.1} There are two vertices $u,w\in V\setminus (F_{1}\bigcup F_{2})$ and there is a vertex $v\in F_{1}\Delta F_{2}$ such that $uw\in E$ and $vw\in E$.

\item [{\rm (2)}]\label{4.2} There are two vertices $u,v\in F_{1}\setminus F_{2}$ and there is a vertex $w\in V\setminus (F_{1}\bigcup F_{2})$ such that $uw\in E$ and $vw\in E$.
\item [{\rm (3)}]\label{4.3} There are two vertices $u,v\in F_{2}\setminus F_{1}$ and there is a vertex $w\in V\setminus (F_{1}\bigcup F_{2})$ such that $uw\in E$ and $vw\in E$.
\end{enumerate}
\end{lem}

\begin{theorem}\label{thm4}
Let $HCN_{n}$ be the $n$-dimensional hierarchical cubic network for $n\geq 2$, then the $g$-good neighbour diagnosability of $HCN_{n}$ under the $MM^{*}$ model satisfies $t_{g}(HCN_{n})\leq 2^{g}(n+2-g)-1$ for $n\geq 2$ and $1\leq g\leq n-1$.
\end{theorem}

\f {\bf Proof.} Let $x_{1}Q_{n}, x_{2}Q_{n}, \cdots, x_{2^{n}}Q_{n}$ be the $2^{n}$ $n$-cubes of $HCN_{n}$. Let $X\subseteq V(x_{1}Q_{n})$ and $HCN_{n}[X]\cong Q_{g}$. Let $F_{1}=N_{HCN_{n}}(X)$ and $F_{2}=N_{HCN_{n}}[X]$. By Lemma~\ref{lem4}, both $F_{1}$ and $F_{2}$ are $g$-good neighbour faulty sets. As $|N_{HCN_{n}}(X)|=2^{g}(n+1-g)$, we have
\vskip0.1cm
~~~~~~$|F_{1}|=2^{g}(n+1-g)\leq 2^{g}(n+2-g)$
\vskip0.1cm
and
\vskip0.1cm
~~~~~~$|F_{2}|=|F_{1}|+|X|$
\vskip0.1cm
~~~~~~~~~~~~$=2^{g}(n+1-g)+ 2^{g}$
\vskip0.1cm
~~~~~~~~~~~~$\leq 2^{g}(n+2-g)$
\vskip0.1cm
As $F_{1}\Delta F_{2}=X$, there is no cross edge between $V(HCN_{n})-(F_{1}\cup F_{2})$ and $F_{1}\Delta F_{2}$. By Lemma~\ref{lem4.1}, the $g$-good neighbour faulty sets of $F_{1}$ and $F_{2}$ are indistinguishable. By Definition~\ref{defi2}, then hierarchical cubic network $HCN_{n}$ is not $g$-good neighbour $2^{g}(n+2-g)$-diagnosable under the $MM^{*}$ model. Thus, $t_{g}(HCN_{n})\leq2^{g}(n+2-g)-1$.
\hfill\qed

\begin{theorem}\label{thm5}
Let $HCN_{n}$ be the $n$-dimensional hierarchical cubic network for $n\geq 2$, then the $g$-good neighbour diagnosability under the $MM^{*}$ model satisfies $t_{g}(HCN_{n})\geq 2^{g}(n+2-g)-1$ for $n\geq 2$ and $1\leq g\leq n-1$.
\end{theorem}

\f {\bf Proof.} Let $x_{1}Q_{n}, x_{2}Q_{n}, \cdots, x_{2^{n}}Q_{n}$ be the $2^{n}$ $n$-cubes of $HCN_{n}$. To Prove the result, we just need to show that for any two distinct $g$-good neighbour faulty subsets $F_{1}$ and $F_{2}$ of $HCN_{n}$ such that $|F_{1}|\leq 2^{g}(n+2-g)-1$ and $|F_{2}|\leq 2^{g}(n+2-g)-1$, the sets $F_{1}$ and $F_{2}$ are distinguishable.

We prove the result by contradiction. That is, there are two distinct $g$-good neighbour faulty subsets $F_{1}$ and $F_{2}$ of $HCN_{n}$ with $|F_{1}|\leq 2^{g}(n+2-g)-1$ and $|F_{2}|\leq 2^{g}(n+2-g)-1$, but they are indistinguishable. Without loss of generality, assume that $F_{2}\setminus F_{1}\neq \emptyset$. To obtain a contradiction, the following fact and claim are useful.
\vskip0.1cm
\f {\bf Fact 1.} With a similar proof as Theorem~\ref{thm2}, we obtain that $V(HCN_{n})\neq F_{1}\bigcup F_{2}$.
\vskip0.2cm
{\bf Claim 1}. There is no isolated vertex in $HCN_{n}\setminus (F_{1}\bigcup F_{2})$.
\vskip0.1cm
\f {\bf Proof of Claim 1.} To prove the result, the following two cases are considered.

Case 1. $2\leq g\leq n-1$.

Suppose to the contrary. That is, $HCN_{n}\setminus (F_{1}\bigcup F_{2})$ has at least one isolated vertex, say $u$. Obviously, we have $d_{HCN_{n}\setminus (F_{1}\bigcup F_{2})}(u)=0$ and

$d_{HCN_{n}[F_{2}\setminus F_{1}]}(u)=d_{HCN_{n}[V(HCN_{n})\setminus F_{1}]}(u)$.
\vskip0.1cm
As $F_{1}$ is a $g$-good neighbour faulty set, then  $d_{HCN_{n}[V(HCN_{n})\setminus F_{1}]}(u)\geq g$ and
\vskip0.1cm
$d_{HCN_{n}[F_{2}\setminus F_{1}]}(u)=d_{HCN_{n}[V(HCN_{n})\setminus F_{1}]}(u)\geq g\geq 2$,
\vskip0.1cm
which satisfies condition $(3)$ of Lemma~\ref{lem4.1}. Thus, the $g$-good neighbour faulty sets $F_{1}$ and $F_{2}$ are distinguishable, a contradiction.

Case 1. $g=1$ and $n\geq 2$.

Suppose to the contrary. That is, $HCN_{n}\setminus (F_{1}\bigcup F_{2})$ has at least one isolated vertex, say $w$. Let $W$ be the set of all isolated vertices in $HCN_{n}\setminus (F_{1}\bigcup F_{2})$ and $H=HCN_{n}\setminus (F_{1}\bigcup F_{2}\cup W)$.

If $F_{1}\setminus F_{2}=\emptyset$, then $F_{1}\subseteq F_{2}$. As $F_{2}$ is a $1$-good neighbour faulty set, $d_{HCN_{n}[V(HCN_{n})\setminus F_{2}]}(w)\geq 1$. As $F_{1}\subseteq F_{2}$, then

$d_{HCN_{n}\setminus (F_{1}\cup F_{2})}(w)=d_{HCN_{n}[V(HCN_{n})\setminus F_{2}]}(w)\geq g\geq 1$,

which contradicts with the fact that $w$ is an isolated vertex in $HCN_{n}\setminus (F_{1}\bigcup F_{2})$.

Now, suppose that $F_{1}\setminus F_{2}\neq\emptyset$. Recall that $w$ is an isolated vertex in $HCN_{n}\setminus (F_{1}\bigcup F_{2})$. Obviously, we have $d_{HCN_{n}[V(HCN_{n})\setminus (F_{1}\cup F_{2})]}(w)=0$ and

$d_{HCN_{n}[F_{1}\setminus F_{2}]}(w)=d_{HCN_{n}[V(HCN_{n})\setminus F_{2}]}(w)$.

As $F_{2}$ is a $1$-good neighbour faulty set, we have

$d_{HCN_{n}[F_{1}\setminus F_{2}]}(w)=d_{HCN_{n}[V(HCN_{n})\setminus F_{2}]}(w)\geq g=1$.

If $d_{HCN_{n}[F_{1}\setminus F_{2}]}(w)\geq 2$, then it satisfies condition (2) of Lemma~\ref{lem4.1}. Then the $g$-good neighbour faulty sets $F_{1}$ and $F_{2}$ are distinguishable, a contradiction.

Thus, $d_{HCN_{n}[F_{1}\setminus F_{2}]}(w)=1$. Let $u\in F_{1}\setminus F_{2}$ such that $uw\in E(HCN_{n})$. Similarly, as $F_{1}$ is a $1$-good neighbour faulty set, we have $d_{HCN_{n}[F_{2}\setminus F_{1}]}(w)=1$. Let $v\in F_{2}\setminus F_{1}$ such that $vw\in E(HCN_{n})$. Thus, $d_{HCN_{n}[F_{1}\cap F_{2}]}(w)=(n+1-2)$ and $|F_{1}\cap F_{2}|\geq (n+1)-2$. For $g=1$, $|F_{2}|\leq 2^{g}(n+2-g)-1=2n+1$. Then we deduce that

$\sum_{w\in W}|N_{HCN_{n}[F_{1}\cap F_{2}](w)}|=|W|[(n+1)-2]$

~~~~~~~~~~~~~~~~~~~~~~~~~~~~~~~~~~$\leq \sum_{v\in F_{1}\cap F_{2}}d_{HCN_{n}}(v)$

~~~~~~~~~~~~~~~~~~~~~~~~~~~~~~~~~~$=|F_{1}\cap F_{2}|(n+1)$

~~~~~~~~~~~~~~~~~~~~~~~~~~~~~~~~~~$\leq(|F_{2}|-1)(n+1)$

~~~~~~~~~~~~~~~~~~~~~~~~~~~~~~~~~~$=2n(n+1)$

Thus, we have $|W|\leq 2n(n+1)/(n-1)$.

If $H=\emptyset$, then

~~~~~~~$|V(HCN_{n})|=2^{2n}$

~~~~~~~~~~~~~~~~~~~~~~~$=|F_{1}\bigcup F_{2}|+|W|$

~~~~~~~~~~~~~~~~~~~~~~~$=|F_{1}|+|F_{2}|-|F_{1}\cap F_{2}|+|W|$

~~~~~~~~~~~~~~~~~~~~~~~$\leq 2(2n+1)-(n+1-2)+2n(n+1)/(n-1)$

~~~~~~~~~~~~~~~~~~~~~~~$=3n+3+2n(n+1)/(n-1)$

However, $2^{2n}> 3n+3+2n(n+1)/(n-1)$ for $n\geq 3$, a contradiction. Thus, $H\neq \emptyset$.

For any vertex $b_{1}\in H$, as $H$ has no isolated vertex, then there exists some vertex $b_{2}\in H$ such that $b_{1}b_{2}\in E(HCN_{n})$. If $b_{1}v_{1}\in E(HCN_{n})$ for some vertex $v_{1}\in F_{1}\bigtriangleup F_{2}$, then the $1$-good neighbour faulty sets $F_{1}$ and $F_{2}$ are distinguishable by condition (1) of Lemma~\ref{lem4.1}, which is a contradiction. Thus, $b_{1}v_{1}\notin E(HCN_{n})$ for any vertex $v_{1}\in F_{1}\bigtriangleup F_{2}$.

By the arbitrariness of $v_{1}$ and $b_{1}$, there is no edge between $H$ and $F_{1}\bigtriangleup F_{2}$.

As both $F_{1}$ and $F_{2}$ are $1$-good neighbour faulty sets, $F_{1}\cap F_{2}$ is a $1$-good neighbour faulty set. As there is no edge between $H$ and $F_{1}\bigtriangleup F_{2}$, $F_{1}\cap F_{2}$ is a $1$-good neighbour faulty cut. By Lemma~\ref{lem3}, $|F_{1}\cap F_{2}|\geq 2(n+1-1)=2n$.

As $|F_{1}|\leq 2n+1, |F_{2}|\leq 2n+1, F_{1}\setminus F_{2}\neq\emptyset$ and $F_{2}\setminus F_{1}\neq\emptyset$, we have $|F_{1}\setminus F_{2}|=1$ and $|F_{2}\setminus F_{1}|=1$.

Let $F_{1}\setminus F_{2}=\{u\}$ and $F_{2}\setminus F_{1}=\{v\}$. By Lemma~\ref{lem6}, $u$ and $v$ have at most two common neighbours. Thus, $|W|\leq 2$ as any vertex of $W$ is adjacent to both $u$ and $v$.

If $|W|=2$, let $W=\{w_{1}, w_{2}\}$. By Lemma~\ref{lem6} and Lemma~\ref{lem2.1.1}, no two vertices of $u, v, w_{1}$ and $w_{2}$ have a common neighbour in $F_{1}\cap F_{2}$. Thus,

~~~~~~$|F_{1}\cap F_{2}|\geq |N(u)\setminus\{w_{1},w_{2}\}|+|N(v)\setminus\{w_{1},w_{2}\}|$

~~~~~~~~~~~~~~~~~~$+|N(w_{1})\setminus\{u,v\}|+|N(w_{2})
\setminus\{u,v\}|$

~~~~~~~~~~~~~~~~~~$=4(n+1-2)$

~~~~~~~~~~~~~~~~~~$=4n-4$.

Then we have

~~~~$|F_{2}|=|F_{1}\cap F_{2}|+|F_{2}\setminus F_{1}|$

~~~~~~~~~~~$\geq 4n-4+1$


~~~~~~~~~~~$> 2n+1$

~~~~~~~~~~~$\geq |F_{2}|$ for $n\geq 2$,

which is a contradiction.

If $|W|=1$, say $W=\{w\}$, then we have $uv\notin E(HCN_{n})$ by Lemma~\ref{lem2.1.1}, then we have

~~~~~~$|N(u,v,w)|\geq |N(u)\setminus{w}|+|N(v)\setminus{w}|+
|N(w)\setminus\{u,v\}|$

~~~~~~~~~~~~~~~~~~~~~~$-|(N(u)\cap N(v))\setminus\{w\}|-|N(u)\cap N(w)|$

~~~~~~~~~~~~~~~~~~~~~~$-|N(v)\cap N(w)|+|N(u)\cap N(v)\cap N(w)|$

~~~~~~~~~~~~~~~~~~~~~~$=2n+(n+1-2)-1+0$

~~~~~~~~~~~~~~~~~~~~~~$=3n-2$

Thus, we have $|F_{1}\cap F_{2}|\geq |N(u,v,w)|\geq 3n-2$

and

~~~~$|F_{2}|=|F_{1}\cap F_{2}|+|F_{2}\setminus F_{1}|$

~~~~~~~~~~~$\geq 3n-2+1$


~~~~~~~~~~~$> 2n+1$

~~~~~~~~~~~$\geq |F_{2}|$ for $n\geq 3$,

which is a contradiction.

The proof of Claim 1 is complete.

By Claim $1$ and Claim $2$, for any vertex $u\in HCN_{n}\setminus (F_{1}\cup F_{2})$, there exists some vertex $v\in HCN_{n}\setminus (F_{1}\cup F_{2})$ such that $uv\in E(HCN_{n})$. If $uw\in E(HCN_{n})$ for $w\in F_{1}\bigtriangleup F_{2}$, it satisfies condition $(3)$ of Lemma~\ref{lem4.1}. Thus, the $g$-good neighbour faulty sets $F_{1}$ and $F_{2}$ are distinguishable, a contradiction. That is to say, $uw\notin E(HCN_{n})$. By the arbitrariness of $u,w$, there is no edge between $HCN_{n}\setminus(F_{1}\cup F_{2})$ and $F_{1}\bigtriangleup F_{2}$, we have

$d_{HCN_{n}[V(HCN_{n})\setminus F_{1}]}(w)=d_{HCN_{n}[F_{2}\setminus F_{1}]}(w)\geq g$.

Thus, by Lemma~\ref{lem2.1}, $|F_{2}\setminus F_{1}|\geq 2^{g}$.

As both $F_{1}$ and $F_{2}$ are $g$-good neighbour faulty sets, $F_{1}\cap F_{2}$ is a $g$-good neighbour faulty set. In addition, as there is no edge between $HCN_{n}\setminus(F_{1}\cup F_{2})$ and $F_{1}\bigtriangleup F_{2}$, $F_{1}\cap F_{2}$ is a $g$-good neighbour faulty cut. By lemma 5, $|F_{1}\bigtriangleup F_{2}|\geq 2^{g}(n+1-g)$. Thus,

~~~~$|F_{2}|=|F_{1}\cap F_{2}|+|F_{2}\setminus F_{1}|$

~~~~~~~~~~~$\geq 2^{g}(n+1-g)+2^{g}$

~~~~~~~~~~~$=2^{g}(n+2-g)$

which contradicts with $|F_{2}|< 2^{g}(n+2-g)-1$.

The proof of the theorem is complete.
\hfill\qed
\vskip0.1cm
By Theorem~\ref{thm4} and Theorem~\ref{thm5}, the following theorem can be obtained.

\begin{theorem}\label{thm6}
Let $HCN_{n}$ be the $n$-dimensional hierarchical cubic network for $n\geq 2$, then the $g$-good neighbour diagnosability under the $MM^{\ast}$ model is $t_{g}(HCN_{n})= 2^{g}(n+2-g)-1$ for $n\geq 2$ and $1\leq g\leq n-1$.
\end{theorem}

\section{Concluding remarks}
As the hierarchical cubic network $HCN_{n}$ has some attractive properties to design interconnection networks. In this paper, we focus on the graph $HCN_{n}$. We show that the $g$-good-neighbour diagnosability of $HCN_{n}$ under the $PMC$ model and $MM^{\ast}$ model is $2^{g}(n+2-g)-1$, respectively. In the future work, we would like to study $g$-extra connectivity and the $g$-extra diagnosability of $HCN_{n}$ under the $PMC$ model and $MM^{\ast}$ model for larger $g$. This could be worth of further investigation.

\section*{Acknowledgments}
This work was supported by the National Natural Science Foundation of China (Nos.11371052, 11731002, 11571035), the Fundamental Research Funds for the Central Universities (No. 2016JBM071, 2016JBZ012) and the $111$ Project of China (B16002).


\begin{thebibliography}{10}


\bibitem{chi} W.-K. Chiang, R.-J. Chen, Topological properties of hierarchical cubic networks, J. Syst. Architecture 42 (4) (1996) 289-307.

\bibitem{D} A.T. Dahbura and G.M. Masson, An O$(n^{2.5})$ fault identification algorithm for diagnosable systems, IEEE Trans. Compu. 33(6) (1984) 486-492.

\bibitem{E} A. Esfahanian, Generalized measures of fault tolerance with application to n cube networks, IEEE Trans. Compu. 38 (1989) 1586-1591.

\bibitem{F} J. Fan, Diagnosability of crossed cubes under the comparison diadnosis model, IEEE Trans. Parallel Distrib. Syst. 13(10) (2002) 1099-1104.

\bibitem{fu1} J.-S. Fu, G.-H. Chen, Hamiltonicity of the hierarchical cubic network, Theory Comput. Syst. 35(1) (2002) 59每79.

\bibitem{fu2} J.-S. Fu, G.-H. Chen, D.-R. Duh, Node-disjoint paths and related problems on hierarchical cubic networks, Networks 40 (2002) 142每154.

\bibitem{gho} K. Ghose, K.R. Desai, Hierarchical cubic network, IEEE Trans. Parallel Distrib. Syst. 6 (1995) 427每435.

\bibitem{Li} X. Li, M. Liu, Z. Yan and J. Xu, On conditional fault tolerance of hierarchical cubic networks, Theore. Compu. Sci., in Press.

\bibitem{L} P.-L. Lai, J. Tan, C. Chang and L. Hsu, Conditional diagnosability measures for large multiprocessor systems, IEEE Trans. Compu. 54(2) (2005) 165-175.

\bibitem{La}  S, Latifi, M. Hedge, M. Naraghi-pour, Conditional connectivity measures for large multiprocessor systems, IEEE Trans. Compu. 43 (1994) 218-222.

\bibitem{M} J. Maeng and M. Malek, A comparison connection assignment for self-diagnosis of multiprocessor systems, in Proceedings of the 11-th International Symp. Fault-Tolerant Computing, T.M. Conte, ed. IEEE Computer Society Press, Los Alamitos, 1981,  173-175.

\bibitem{O} A. Oh and H. Chol, Generalized measures of fault tolerance in n-cube networks, IEEE Trans. Parallei Dist. Syst. 4 (1993) 702-703.

\bibitem{Pe} S. Peng, C. Lin, J.M. Tan and L, Hsu, The $g$-good-neighbour conditional diagnosability of hypercube under the PMC model, App. Math. Comput. 218(21) (2012) 10406-10412.

\bibitem{Pr}F. Preparata, G. Metze and R. Chien, On the connection assighment problem of diagnosis systems, IEEE Trans. Comput. EC-16 (1967) 848-854.


\bibitem{Wan} M. Wan and Z. Zhang, A kind of conditional connectivity of star graphs, Appl. Math. Lett. 22 (2009) 264-267.


\bibitem{Wu} J. Wu, G. Guo, Fault tolerance measures for m-ary n-dimensional hypercubes based on forbidden faulty sets, IEEE Trans. Comput. 47 (1998) 888每893.

\bibitem{w} S. Wang and W. Han, The $g$-good-neighbour conditional diagnosability of $n$-dimensional hypercubes under the $MM^{\ast}$ model, Inform. Process. Lett. 116 (2016) 574-577.

\bibitem{w1} M. Wang, Y. Guo and S. Wang, The $1$-good-neighbour diagnosability of Cayley graphs generated by transposition trees under the PMC model and $MM^{\ast}$ model, International J. of Compu. Math. 94 (2017) 620-631.

\bibitem{w2} M. Wang, Y. Lin and S. Wang, The $2$-good-neighbour diagnosability of Cayley graphs generated by transposition trees under the PMC model and $MM^{\ast}$ model, Theoret. Compu. Sci. 628 (2016) 92-100.

\bibitem{w3} M. Wang, Y. Lin and S. Wang, The $1$-good-neighbour connectivity and diagnosability of Cayley graphs generated by complete graphs, Discrete Appl. Math. 295 (2017) 536-543.

\bibitem{X} J. Xu, Topological Structure and Analysis of Interconnection network, Kluwer Academic Pbulishers, Dordrecht/Boston/London. 2001.

\bibitem{Y} J. Yuan, A, Liu X. Ma, X. Qin and J. Zhang, The $g$-good-neighbour conditional diagnosability of $k$-ary $n$-cubes under the $PMC$ and $MM^{*}$ model , IEEE Trans. Parallel Distrib. Syst. 26(10) (2015) 1165-1177.


\bibitem{Yang} W. Yang, H. Li and J. Meng, Conditional connectivity of Cayley graphs generated by transposition trees, Inf. Process. Lett. 110 (2010) 1027-1030.

\bibitem{Yu} X. Yu, X. Huang and Z. Zhang, A kind of conditional connectivity of Cayley graphs generated by unicyclic graphs, Inf. Sci. 243 (2013) 86-94.

\bibitem{Yang1} Z. Zhang, W. Xiong and W. Yang, A kind of conditional fault tolerance of alternating group graphs, Inf. Process. Lett. 110 (2010) 998-1002.

\bibitem{yun1} S.-K. Yun, K.-H. Park, The optimal routing algorithm in hierarchical cubic network and its properties, IEICE Trans. Inform. Syst. 78(4) (1995) 436-443.

\bibitem{yun2} S.-K. Yun, K.-H. Park, Comments on hierarchical cubic network, IEEE Trans. Parallel Distrib. Syst. 9 (1998) 410-414.

\bibitem{zhou} S. Zhou, S. Song, X. Yang, L. Chen, On conditional fault tolerance and diagnosability of hierarchical cubic networks, Theoret. Comput. Sci. 609 (2016) 421-433.

\bibitem{Zhu} Q. Zhu, J.M. Xu, On restricted edge connectivity and extra edge connectivity of hypercubes and folded hypercubes, J. Univ. Sci. Technol. China 36 (3) (2006) 246-253.

\end{thebibliography}
\end{document}